\documentclass[a4paper,12pt]{amsart} 

\usepackage{amsmath,graphics}
\usepackage{amssymb}
\usepackage{amsfonts}
\usepackage{latexsym}
\usepackage{eucal}
\usepackage[dvips]{graphicx}
\usepackage{mathrsfs}
\usepackage{enumerate}
\newtheorem{thm}{Theorem}[section]

\newtheorem{propo}[thm]{Proposition}
\newtheorem{lem}[thm]{Lemma}
\newtheorem{cor}[thm]{Corollary}

\newtheorem{conj}[thm]{Conjecture}
\renewcommand{\Re}{{\rm Re}}
\renewcommand{\Im}{{\rm Im}}
\newcommand{\R}{\mathbb{R}}
\newcommand{\C}{\mathbb{C}}
\newcommand{\F}{\mathbb{F}}

\newcommand{\Z}{\mathbb{Z}}
\newcommand{\N}{\mathbb{N}}

\renewcommand{\H}{\mathbb{H}^2}

\newcommand{\lt}{{\mathcal L}}

\newcommand{\D}{{\mathcal D}}
\newcommand{\G}{{\mathbf G}}

\newcommand{\half}{{\textstyle{\frac{1}{2}}}}
\newcommand{\halfdelta}{{\textstyle{\frac{\delta}{2}}}} 
\begin{document}
\bibliographystyle{plain}
\title[Resonances and congruence subgroups]{ Resonances and convex co-compact congruence subgroups of $PSL_2(\Z)$}

\author[D. Jakobson]{Dmitry Jakobson}
\address{McGill University\\
Department of Mathematics and Statistics\\
805 Sherbrooke Street West\\
Montreal, Quebec, Canada H3A0B9 
}
\email{jakobson@math.mcgill.ca}

\author[F.~Naud]{Fr\'ed\'eric Naud}
\address{%
Fr\'ed\'eric Naud\\
Laboratoire d'Analyse non-lin\'eaire et G\'eom\'etrie\\
Universit\'e d'Avignon, 33 rue Louis Pasteur\\
84000 Avignon\\
France.
}
\email{frederic.naud@univ-avignon.fr}

 \maketitle

 \begin{abstract}
  Let $\Gamma$ be a convex co-compact subgroup of $PSL_2(\Z)$, and let $\Gamma(q)$ be the sequence of "congruence" subgroups of $\Gamma$.
  Let $\mathcal{R}_q\subset \C$ be the resonances of the hyperbolic Laplacian on the "congruence" surfaces $\Gamma(q)\backslash \H$. We prove two results on the  density of resonances in $\mathcal{R}_q$ as $q \rightarrow \infty$: the first shows at least $Cq^3$ resonances in slowly growing discs, the other one is a bound from above in boxes $\{ \delta/2<\sigma\leq \Re(s) \leq \delta\}$, with $\vert \Im(s) -T \vert \leq 1$, where we prove a density estimate of the type $O(T^{\delta-\epsilon_1(\sigma)}q^{3-\epsilon_2(\sigma)})$ with $\epsilon_j(\sigma)>0$ for all $\sigma>\delta/2$, $j=1,2$. \end{abstract}
 
 \tableofcontents

 \section{Introduction and results} 
 Recently, "thin" subgroups of $SL_2(\Z)$ have attracted some attention in Number Theory. By "thin" we mean an {\it infinite index} subgroup $\Gamma \subset SL_2(\Z)$ whose dimension of the limit set $\delta_\Gamma$ satisfies
 $0<\delta_\Gamma<1$. While the results of Bourgain-Gamburd-Sarnak \cite{BGS} focus on the density of almost primes found among entries of the orbits of thin subgroups, the works of Bourgain and Kontorovich \cite{BK1,BK2} are concerned
 with the density of various subsets of $\N$ obtained through the action of $\Gamma$. One of the key steps of the proofs involves reduction (localization) modulo $q$ where $q$ is a square-free integer, in particular one is led to
 consider "congruence subgroups" $\Gamma(q)$ defined by
 $$\Gamma(q):=\left \{ \gamma \in \Gamma\ :\ \gamma \equiv \left ( \begin{array}{cc}1&0\\0&1 \end{array} \right ) \mathrm{mod}\ q  \right \}.$$
 A critical ingredient is the spectral theory of the {\it infinite area} hyperbolic surfaces
  $$X_q:=\Gamma(q) \backslash \H$$
   where uniform estimates on the spectrum of the Laplacian are often required.  Let us recall some basic known facts about the Laplacian on these objects. 
   Let $\H$ be the hyperbolic plane endowed with its standard metric of constant gaussian curvature $-1$.
Let $\Gamma$ be a convex co-compact discrete subgroup of isometries acting on $\H$. This means
that $\Gamma$ admits a finite sided polygonal fundamental domain in $\H$, with infinite area. We will require that $\Gamma$ has no {\it elliptic} elements different from the identity and that $\Gamma$ has no parabolic elements (no cusps). Under these assumptions, the quotient space 
$X=\Gamma \backslash \H$ is a hyperbolic surface with infinite area whose ends are given by hyperbolic funnels. The limit set of $\Gamma$
is commonly defined by 
$$\Lambda(\Gamma):=\overline{\Gamma.z}\cap \partial \H,$$
where $z\in \H$ is a given point and $\Gamma.z$ is the orbit of that point under the action of $\Gamma$ which by discreteness accumulates
on the boundary $\partial \H$. The limit set $\Lambda$ does not depend on the choice of $z$ and its Hausdorff dimension is commonly denoted by $\delta_\Gamma$.
It is the critical exponent of Poincar\'e series \cite{Patterson1} (here $d$ denotes hyperbolic distance)
$$P_\Gamma(s):=\sum_{\gamma \in \Gamma} e^{-s {\mathrm d}(\gamma z,z')}.$$
Let $\Delta$ be the hyperbolic Laplacian on $X$. Its spectrum on $L^2(X)$ has been described completely by 
Lax and Phillips in \cite{LP2}. The half line $[1/4, +\infty)$ is the continuous spectrum and it contains no embedded eigenvalues.
The rest of the spectrum (point spectrum) is empty if $\delta\leq \half$, finite and starting at $\delta(1-\delta)$ if $\delta>\half$. The fact that the bottom of the spectrum is related to the dimension $\delta$ was first discovered by Patterson \cite{Patterson1} for convex co-compact groups.

\bigskip \noindent
By the preceding description of the spectrum, the resolvent 
$$R_\Gamma(s)=\left(\Delta-s(1-s) \right)^{-1}:L^2(X)\rightarrow L^2(X),$$
is therefore well defined and analytic on the half-plane $\{ \Re(s)> \half \}$ except at a possible finite set
of poles corresponding to the finite point spectrum. {\it Resonances} are then defined as poles of the meromorphic continuation of 
$$R_\Gamma(s):C_0^\infty(X)\rightarrow C^{\infty}(X)$$ 
to the whole complex plane. The set of poles is denoted by ${\mathcal R}_X$. This continuation is usually performed via the analytic Fredholm theorem
after the construction of an adequate parametrix. The first result of this kind in the more general setting
of asymptotically hyperbolic manifolds is due to Mazzeo and Melrose \cite{MazzMel}.
A more precise parametrix for surfaces was constructed by Guillop\'e and Zworski in \cite{GuiZwor,GuiZwor1}.
It should be mentioned at this point that in the infinite area case, resonances are spread all over the half plane $\{ \Re(s)< \delta\}$,
in sharp contrast with the finite area non-compact case where resonances are known to be confined in a strip. Among the known results (and conjectures) on the density and location of resonances we mention the following
two facts which are relevant for this paper. Let $\mathcal{N}_X(r)$ be the counting function defined by 
$$\mathcal{N}_X(r):=\# \{ s \in \mathcal{R}_X\ :\  \vert s\vert \leq r \}.$$
From the work of Guillop\'e and Zworski, we know that there exists $C_X>0$ such that for all $r\geq 1$, we have
$$C_X^{-1}r^2\leq \mathcal{N}_X(r)\leq C_X r^2.$$
On the other hand, let $\mathcal{M}_X(\sigma,T)$ be defined by
$$\mathcal{M}_X(\sigma,T):=\# \{ s \in \mathcal{R}_X\ :\ \sigma \leq \Re(s) \leq \delta\ \mathrm{and}\ \vert \Im(s)-T\vert \leq 1 \}.$$
From the work of Guillop\'e-Lin-Zworski \cite{GLZ},  we know that one can find $C_\sigma>0$ such that for all $T \geq 1$, we have
$$ \mathcal{M}_X(\sigma,T)\leq C_\sigma T^{\delta}.$$
It is conjectured in Jakobson-Naud \cite{JakobsonNaud2} that for all $\sigma>\delta/2$ and all $T$ large enough, $\mathcal{M}_X(\sigma,T)=o(1)$, in other words that there exists an "essential" spectral gap up to 
$\{ \Re(s)=\delta/2 \}$ which plays the role of the critical line in infinite volume.

\bigskip \noindent
In this paper, motivated by number theoretic works, we will restrict ourselves to the interesting case when $\Gamma$ is a convex co-compact subgroup of $PSL_2(\Z)$, and will assume throughout the paper that $\Gamma$ is {\it non-elementary} i.e. is not generated by a single hyperbolic element. I should be noticed that since $\Gamma$ is a {\it free} group, there is no need to distinguish
$\Gamma$ as a subgroup of $PSL_2(\Z)$ or viewed as a matrix subgroup of $SL_2(\Z)$. \footnote{Indeed, if $\Gamma$ is a free subgroup 
of $SL_2(\Z)$, $-Id\not \in \Gamma$ and therefore the natural projection $SL_2(\Z)\rightarrow PSL_2(\Z)$ is injective when restricted to $\Gamma$.}
As mentioned above, a natural question is to describe the resonances of "congruence" surfaces $\Gamma(q)$ with respect to $q$. For simplicity, we will restrict ourselves in this paper to the case when {\it $q$ is a prime number}. From the work of Gamburd \cite{Gamburd}, it is known that for
all $q$ large enough, the projection map
$$\Phi_q:\left \{  \Gamma \rightarrow SL_2(\F_q) \atop \gamma \mapsto \gamma\  \mathrm{mod}\  q \right. ,$$
is a surjection. Therefore we have
$$[\Gamma:\Gamma(q)]=\vert SL_2(\F_q) \vert=q(q^2-1)\asymp q^3$$ 
Since all the subgroups $\Gamma(q)$ have a finite index in $\Gamma$, they have all the same dimension $\delta_{\Gamma(q)}=\delta_\Gamma$.
Our first result is the following. For simplicity, we denote the counting function $\mathcal{N}_{\Gamma(q)}(r)$ by $\mathcal{N}_q(r)$.
\begin{thm}
\label{lowerb}
There exist constants $C_0>0$ and $T_0>0$ such that for all $\epsilon>0$ and all $q\geq q_0(\epsilon)$ and $T\geq T_0$, we have 
$$\mathcal{N}_q\left(T(\log q)^\epsilon \right)\geq C_0 T^2 q^3.$$
\end{thm}
This result shows abundance of resonances in discs with slow radius growth as $q\rightarrow \infty$. This lower bound is not surprising in view of the geometric bounds obtained by Borthwick in \cite{Borthwick3}. Notice that this bound is "almost optimal" in term of $q$, indeed, in $\S 2$, we show that
$$ \mathcal{N}_q(T)\leq Cq^3 \log(q) T^2,$$
uniformly for all $T\geq 1$.
It may be possible (with some more geometric work) to derive Theorem \ref{lowerb} directly from the lower bound in \cite{Borthwick3}. However, in this paper, we use a fairly different (and more algebraic) route which is justified by our next statement.
We denote by $\mathcal{R}_q$ the resonance spectrum of the surface $X_q=\Gamma(q)\backslash \H$.

\noindent In their work on almost primes \cite{BGS}, Bourgain-Gamburd-Sarnak (and also Gamburd \cite{Gamburd}) obtained the following "spectral gap" result for the family $\Gamma(q)$:

\begin{itemize}
\item If $\delta>1/2$ then there exists $\epsilon_0>0$, independent of $q$ such that for all $q$  we have $\mathcal{R}_q\cap \{ \Re(s) >\delta-\epsilon_0\}=\{ \delta\}$.
\item Moreover, if $\delta>5/6$ then $\mathcal{R}_q\cap [5/6, \delta]$ is independent of $q$ (notice that in this case there are only eigenvalues).
\item If $\delta\leq 1/2$, then there exists $\epsilon_0>0$, independent of $q$ such that for all $q$ we have
$$ \mathcal{R}_q\cap \left \{ \Re(s)>\delta-\epsilon_0\min \left(1, \frac{1}{\log(1+\vert \Im(s) \vert)} \right ) \right \}=\{ \delta \}.$$
\end{itemize}
These uniform spectral gaps are obtained thanks to the work of Bourgain-Gamburd \cite{BG1} on Cayley graphs of $SL_2(\F_q)$ which are proved to be expanders for all finitely generated, non elementary set of generators.
While in the case $\delta>1/2$ they can use the full strength of Lax-Phillips spectral theory, in the case $\delta\leq 1/2$, one has to face genuine resonances and use some transfer operators techniques from \cite{Naud2}. In view of the conjecture of \cite{JakobsonNaud2}, it is natural to expect the following {\it "uniform essential spectral gap property"}.

\begin{conj}
For all $\sigma>\halfdelta$, $\mathcal{R}_q\cap \{ \Re(s)\geq \sigma \}$ is finite and independent of $q$.
\end{conj}
This conjecture seems to be out of reach (in the finite volume case, this is Selberg's eigenvalue conjecture), but we can prove the following. 

\begin{thm}
 \label{upperb}
  Assume that $\sigma>\halfdelta$, then there exist $m_0(\sigma)$, such that we have for all $T\geq 1$ and $q$ large,
  $$ \mathcal{M}_q(\sigma,T):=\# \{ s \in \mathcal{R}_q\ :\ \sigma \leq \Re(s) \leq \delta\ \mathrm{and}\ \vert \Im(s)-T \vert \leq 1 \}$$
  $$\leq m_0 T^{\delta+\tau_1(\sigma)} q^{3+\tau_2(\sigma)},$$
  where for $i=1,2$, $\tau_i(\sigma)<0$ on $(\halfdelta,\delta]$, and $\tau_i$ is strictly convex and decreasing on $(\halfdelta, \delta]$.
\end{thm}
\noindent This statement is a strengthening of the main result in \cite{NaudInvent}, extended to all congruence subgroups. Not only we have a gain over the crude bound
$O(T^{\delta})$ but simultaneously a gain over the $O(q^3)$ bound, as long as we count resonances in $\{ \Re(s)\geq \sigma >\halfdelta \}$.
Notice that it can be rephrased (for fixed $T$) as a bound of the type 
 $$O\left( \mathrm{Vol}(N_q)^{1+\frac{\tau_2(\sigma)}{3}} \right),$$
 where $N_q$ is the Nielsen region in $\Gamma(q) \backslash \H$ (the convex core of the surface).
 
In the particular case $\delta>\half$ it gives a density theorem for the number of eigenvalues of the Laplacian $\Delta_q$ in $[\delta(1-\delta),1/4]$ which to our knowledge does not follow from previous works:

$$ \# (\mathrm{Sp}(\Delta_q)\cap [\delta(1-\delta),1/4]) =O\left ( \mathrm{Vol}(N_q)^{1-\epsilon_\Gamma} \right )$$
for some $\epsilon_\Gamma>0$ depending only on $\Gamma$. 
 This resonance behavior in the strip $\{ \delta/2< \Re(s)\leq \delta\}$ is in sharp contrast with Theorem \ref{lowerb}, 
which shows a drastically different behavior in the half-plane $\{ \Re(s)\leq \delta/2 \}$. We point out that the functions $\tau_i(\sigma)$ have an "explicit" formula in terms of topological pressure. We do not expect this formula to produce a uniform spectral gap for $\sigma$ close to 
$\delta$, although we use some of the ideas of \cite{BGS, Gamburd} in our proof. 

Let us describe the organization of the paper. The main tool in both results is to use "congruence" transfer operators which were already defined in \cite{BGS}. In section $\S2$ we recall how  convex co-compact subgroups of $PSL_2(\Z)$ can be viewed as Schottky groups. We then define the so-called congruence transfer operator on a suitable space of Holomorphic functions and 
show that its Fredholm determinant is related to the Selberg's zeta function of $X_q$. An upper bound on the growth of this determinant is then proved, using some singular values estimates.
Combining this result with the Trace formula leads to Theorem \ref{lowerb}, see $\S 3$. To prove Theorem \ref{upperb}, a different kind of approach is required. 
In $\S 4$, we use some (modifications of) ideas from \cite{NaudInvent}, where we estimate the number of resonances by using a regularized Hilbert-Schmidt determinant related to iterates of the transfer operator. A key part of the proof comes from a variant on the lower bound on the girth of Cayley graphs of $SL_2(\F_q)$ which is proved in \cite{Gamburd} and a separation mechanism based on the disconnected topology of the limit set. 

\bigskip 
\noindent {\bf Acknowledgments.} This work was completed while FN was visiting CRM at universit\' e de Montr\'eal under a CNRS funding. Both authors are supported by ANR "blanc" GeRaSic.
DJ is also supported by NSERC, FQRNT and Peter Redpath Fellowship.

\section{Congruence transfer operators and Selberg's zeta functions}

\subsection{The Schottky picture}
We use the notations of $\S 1$. Let $\H$ denote the Poincar\'e upper half-plane $\H=\{ x+iy\in \C\ :\ y>0\}$ endowed with its standard metric of constant curvature $-1$
$$ds^2=\frac{dx^2+dy^2}{y^2}.$$ 
The group of isometries of $\H$ is $\mathrm{PSL}_2(\R)$ through the action of 
$2\times 2$ matrices viewed as M\"obius transforms
$$z\mapsto \frac{az+b}{cz+d},\ ad-bc=1.$$ 
Below we recall the definition of Fuchsian Schottky groups which will be used to define transfer operators.
A Fuchsian Schottky group is a free subgroup of $\mathrm{PSL}_2(\R)$ built as follows. Let $\D_1,\ldots, \D_p,\D_{p+1},\ldots, \D_{2p}$
be $2p$ Euclidean {\it open} discs in $\C$ orthogonal to the line $\R\simeq \partial \H$. We assume that for all $i\neq j$, $\overline{\D_i} \cap \overline{\D_j}=\emptyset$. 
Let $\gamma_1,\ldots,\gamma_p \in \mathrm{PSL}_2(\R)$ be $p$ isometries such that for all $i=1,\ldots,p$, we have
$$\gamma_i(\D_i)=\widehat{\C}\setminus \overline{\D_{p+i}},$$
where $\widehat{\C}:=\C\cup \{ \infty \}$ stands for the Riemann sphere.

\bigskip




\bigskip \noindent
Let $\Gamma$ be the free group generated by $\gamma_i,\gamma_i^{-1}$ for $i=1,\ldots,p$, then $\Gamma$ is a convex co-compact group, i.e. it is finitely generated
and has no non-trivial parabolic element. The {\it converse is true}: up to an isometry, all convex co-compact hyperbolic surfaces
can be uniformized by a group as above, see \cite{button}. In the particular case when $\Gamma$ is a convex co-compact subgroup of $PSL_2(\Z)$, then by using the same argument as in \cite{button}, one can find a set of generators as above.

\subsection{Topological pressure and Bowen's formula}
For all $j$
we set $I_j:=\D_j\cap \R$ and define the Bowen-Series map $B:\cup_{j=1}^{2p} I_j \rightarrow \R\cup \{\infty\}$ by 
$$B(x):=\gamma_j(x)\ \mathrm{if}\ x\in I_j.$$ 
The maximal $B$-invariant compact subset of $\cup_{j=1}^{2p} I_j$ is precisely the limit set $\Lambda(\Gamma)$, and $B$ is uniformly expanding on $\Lambda(\Gamma)$.
The {\it topological pressure} $P(x)$, $x\in \R$,  is the thermodynamical quantity given by the limit (the sums runs over $n$-periodic points of the map $B$)
\begin{equation}
\label{pressure}
e^{P(x)}=\lim_{n\rightarrow \infty} \left (  \sum_{B^nw=w} \vert (B^n)'(w) \vert^{-x} \right)^{1/n}.
\end{equation}
The fact that this limit exists and defines a {\it real-analytic decreasing strictly convex} function $x\mapsto P(x)$ follows from classical thermodynamical formalism, see for example \cite{PP} for a basic reference, see also \cite{NaudInvent}, for a justification of strict convexity based on the fact that $\Gamma$ is non elementary. Moreover, it has a unique zero on the real line which is exactly the dimension $\delta(\Gamma)$, this is a celebrated result of Bowen \cite{Bowen1}. In particular, we have
$P(x)<0$ iff $x>\delta$, which is something to keep in mind in the rest of the paper, especially in the last section.

\subsection{Determinants and Selberg's zeta functions} 
In the sequel, we will denote by $\G$ the group $SL_2(\F_q)$. Let $\Gamma$ be a convex co-compact subgroup of $PSL_2(\Z)$ as above and
let $\gamma_1,\ldots,\gamma_p$ be a set of Schottky generators as above. For simplicity we denote by $\gamma_{p+i}:=\gamma_i^{-1}$ for $i=1,\ldots,p$.
For each map $\gamma_i$, we fix a $2\times 2$ matrix representation in $SL_2(\Z)$ also denoted by $\gamma_i$. Let 
$$F:\cup_{i=1}^{2p} \D_i\times \G\rightarrow \C$$
be a $\C$-valued function, then the {\it congruence transfer operator} applied to $F$ is defined for all $z\in \D_i,\ g\in \G$
$$\lt_s(F)(z,g):=\sum_{j\neq i} (\gamma_j'(z))^s F(\gamma_j z, \Phi_q(\gamma_j)g ),$$
where $\Phi_q:SL_2(\Z)\rightarrow \G$ is the reduction mod $q$. For obvious simplicity, we will omit $\Phi_q$ in the notations for the right factor.
We need some additional notations. Considering a finite sequence $\alpha$ with
\[\alpha=(\alpha_1,\ldots,\alpha_n)\in \{1,\ldots, 2p\}^n,\]
we set 
\[ \gamma_\alpha:=\gamma_{\alpha_1}\circ \ldots \circ \gamma_{\alpha_n}. \]
We then denote by $\mathscr{W}_n$ the set of admissible sequences of length $n$ by
\[ \mathscr{W}_n:=\left \{ \alpha \in \{1,\ldots, 2p\}^n\ :\ 
\forall\ i=1,\ldots,n-1,\ \alpha_{i+1}\neq \alpha_i +p\ \mathrm{mod}\ 2p \right \}.\]
The set $\mathscr{W}_n$ is simply the set of reduced words of length $n$.
For all $j=1,\ldots, 2p$, we define $\mathscr{W}_n^j$ by
\[ \mathscr{W}_n^j:=\{ \alpha \in \mathscr{W}_n\ :\ \alpha_n\neq j \}. \] 
If $\alpha \in \mathscr{W}_n^j$, then $\gamma_\alpha$ maps $\overline{\D_j}$ into $\D_{\alpha_1+p}$. Using this set of notations, we have the formula for 
$z\in \D_j$,
$$\lt_s^N(F)(z,g)=\sum_{\alpha \in \mathscr{W}_N^j} (\gamma_\alpha'(z))^s F(\gamma_\alpha z, \gamma_\alpha g).$$
We now have to specify the function space on which the transfer operators $\lt_s$ will act. Let $H^2_q$ denote the vector space of (complex-valued) functions $F$
on $\cup_{i=1}^{2p} \D_i\times \G$ such that for all $g\in \G$, $z\mapsto F(z,g)$ is holomorphic on each disc $\D_j$ and such that the following norm
$$ \Vert F \Vert_{q}^2:=\sum_{g \in \G} \int_{\cup_{i=1}^{2p} \D_i} \vert F(z,g)\vert^2dm(z),$$
is finite ($dm$ stands for Lebesgue measure on $\C$). This function space may be viewed as a (vector valued) variant of the classical Bergman spaces, and is a natural Hilbert space. Since each branch $\gamma_i$ acts by contraction on $\cup_{j\neq i} \D_j$, the transfer operators are compact, {\it trace class operators}. This fact is well known and dates back to Ruelle \cite{Ruelle0}, see also Bandtlow-Jenkinson \cite{Bandtlow1} for an in-depth analysis of 
spectral properties of transfer operators on Holomorphic function spaces.
Before we carry on our analysis, it is necessary to recall  a few basic facts on transfer operators acting on $H^2_q$. We start by some distortion estimates.
\begin{itemize}
 \item {\it(Uniform hyperbolicity)}. One can find $C>0$ and $0<\overline{\theta}<\theta<1$ such that for all $n,j$ and $\alpha \in \mathscr{W}_n^j$, for all $z\in \D_j$ we have
 \begin{equation}
 \label{hyp1}
 C^{-1}\overline{\theta}^n\leq \vert \gamma'_\alpha(z) \vert \leq C \theta^n. 
 \end{equation}
 \item {\it (Bounded distortion).} There exists $M_1>0$ such that for all $n,j$ and all $\alpha \in \mathscr{W}_n^j$, for all $z_1,z_2 \in \D_j$ 
 \begin{equation}
\label{bdist}
e^{-\vert z_1-z_2\vert M_1}\leq \frac{\vert \gamma'_\alpha(z_1) \vert}{\vert \gamma'_\alpha(z_2)\vert}\leq  e^{\vert z_1-z_2\vert M_1}.
\end{equation}
\end{itemize}
We refer the reader to \cite{NaudInvent} for details on proofs and references. We will also need the following fact which is proved in \cite{NaudInvent}.
\begin{lem} For all $\sigma_0,M$ in $\R$ with $0\leq \sigma_0<M$, one can find $C_0>0$ such that for all $n$ large enough and $M\geq \sigma\geq \sigma_0$, we have
\begin{equation}
\label{pressure2}
 \sum_{j=1}^{2p}\left (\sum_{\alpha \in \mathscr{W}_n^j} \sup_{I_j}   (\gamma'_\alpha  )^\sigma \right)\leq C_0 e^{nP(\sigma_0)}.
\end{equation}
\end{lem}
With these preliminaries in hand one can prove the following estimate.
\begin{propo}
\label{normest1} There exist a constant $C>0$, independent of $q$ such that for all $N\in \N$, we have 
$$ \Vert \lt_s^N \Vert_{H^2_q}\leq C e^{C\vert s \vert} e^{NP(\Re(s))}.$$
\end{propo}
A straightforward and important consequence is that the {\it spectral radius} of $\lt_s:H^2_q \rightarrow H^2_q$ is bounded by $e^{P(\Re(s))}$. We postpone the proof of this Proposition to the appendix and move on to the central subject of this $\S$.

We recall that the Selberg zeta function $Z_\Gamma(s)$ is defined as the analytic continuation to $\C$ of the infinite product:
$$Z_\Gamma(s):=\prod_{k \in \N} \prod_{\mathcal{C} \in \mathcal{P}} \left (1-e^{-(s+k)\ell(\mathcal{C})} \right),\ \Re(s)>\delta$$
where $\mathcal{P}$ is the set of prime closed geodesics on $\Gamma \backslash \H$, and if $\mathcal{C} \in \mathcal{P}$, $\ell(\mathcal{C})$ is the length.
Our first observation is the following.
\begin{propo}
 \label{detformula} Using the above notation, we have for all $s \in \C$ and $q\geq 2$, 
 $$\det(I-\lt_s)=Z_{\Gamma(q)}(s).$$
\end{propo}

\noindent {\it Proof}. We prove this identity by analytic continuation. First some trace computations are required. 
let $\mathscr{D}_g$ denote the dirac mass at $g\in \G$, i.e.
$$\mathscr{D}_g(h)=\left \{ 1\ \mathrm{if}\  h=g \atop 0\ \mathrm{elsewhere.} \right. $$ 
For all $j=1,\ldots,2p$, we write $\D_j:=D(c_j,r_j)$
and we denote by $\mathbf{e}_k^\ell$ the function  defined for $z\in \D_j$ by
$$\mathbf{e}_k^\ell(z)=\left \{ \begin{array}{c}
0\ \mathrm{if}\ j\neq \ell \\
\sqrt{\frac{k+1}{\pi}}\frac{1}{r_j} \left (\frac{z-c_j}{r_j} \right)^k\ \mathrm{if}\ j=\ell.
\end{array}\right. $$
It is easy to check that the family 
$$\left ( \mathbf{e}_k^\ell  \otimes \mathscr{D}_g \right )_{k\in \N,\ g\in \G}^{\ell=1,\ldots,2p}$$
is a Hilbert basis of $H^2_q$. Writing
$$\mathrm{Tr}(\lt_s^N)=\sum_{k,\ell,g} \langle \lt_s^N(\mathbf{e}_k^\ell  \otimes \mathscr{D}_g), \mathbf{e}_k^\ell  \otimes \mathscr{D}_g \rangle_{H^2_q},$$
we obtain after several applications of Fubini
$$\mathrm{Tr}(\lt_s^N)=\sum_{k,\ell,g} \sum_{\alpha \in \mathscr{W}_N^\ell \atop \alpha_1=p+\ell} \mathscr{D}_g(\gamma_\alpha g) \int_{\D_\ell} (\gamma_\alpha'(z))^s
 \mathbf{e}_k^\ell(\gamma_\alpha z) \overline{\mathbf{e}_k^\ell(z)}dm(z),$$
 $$=\vert \G \vert \sum_{\ell=1}^{2p} \sum_{\alpha \in \mathscr{W}_N^\ell,\alpha_1=p+\ell \atop \gamma_\alpha \equiv \mathrm{Id\ mod}\ q} \sum_{k \in \N}\int_{\D_\ell} (\gamma_\alpha'(z))^s
 \mathbf{e}_k^\ell(\gamma_\alpha z) \overline{\mathbf{e}_k^\ell(z)}dm(z).$$
 We recall that by the mapping property of Schottky groups, there exists $\epsilon_0>0$ such that for all $N$ and all $\alpha \in \mathscr{W}_N^\ell$,
 $$\mathrm{dist}(\gamma_\alpha(\D_\ell),\partial \D_{p+\alpha_1})\geq \epsilon_0.$$
 This uniform contraction property guarantees uniform convergence of 
 $$\sum_k  \mathbf{e}_k^\ell(\gamma_\alpha z) \overline{ \mathbf{e}_k^\ell(z)}$$
 on $\D_\ell \times \D_\ell$, allowing us to write 
 $$ \sum_{k \in \N}\int_{\D_\ell} (\gamma_\alpha'(z))^s
 \mathbf{e}_k^\ell(\gamma_\alpha z) \overline{\mathbf{e}_k^\ell(z)}dm(z)=\int_{\D_\ell} (\gamma_\alpha'(z))^s
 B_{\D_\ell}(\gamma_\alpha z,z) dm(z),$$
where $B_{\D_\ell}(w,z)$ is the Bergman reproducing kernel of the disc $\D_\ell$, given by the explicit formula
$$B_{\D_\ell}(w,z)=\frac{r_\ell^2}{\pi \left [ r_\ell^2-(w-c_\ell)(\overline{z}-c_\ell)\right ]^2}.$$
A standard computation involving Stokes's and Cauchy formula (see  for example Borthwick \cite{Borthwick}, P. 306) then shows that 
$$\int_{\D_\ell} (\gamma_\alpha'(z))^s
 B_{\D_\ell}(\gamma_\alpha z,z) dm(z)= \frac{(\gamma_\alpha'(x_\alpha))^s}{1-\gamma_\alpha'(x_\alpha)},$$
 where $x_\alpha$ is the unique fixed point of $\gamma_\alpha:\D_\ell\rightarrow \D_\ell$. Moreover, 
 $$\gamma_\alpha'(x_\alpha)=e^{-\ell( \mathcal{C}_\alpha)},$$
 where $\mathcal{C}_\alpha$ is the closed geodesic represented by the conjugacy class of $\gamma_\alpha \in \Gamma$. 
 There is a one-to-one correspondence between prime periodic orbits of the Bowen-Series map $B$ and prime conjugacy classes in $\Gamma$ (see Borthwick \cite{Borthwick}, P. 303), therefore
 each prime conjugacy class in $\Gamma$ (and iterates) appears in the above sum over all periodic orbits of $B$. However, for all $\gamma \in \Gamma(q)$ its conjugacy class in $\Gamma$ 
 {\it splits} into 
 $$[\Gamma(q):\Gamma]=\vert G \vert$$
  conjugacy classes in $\Gamma(q)$, with same geodesic length. Let us explain this fact.
 Let $H$ be a normal subgroup of a group $G$, and let $x \in H$. Then it is a basic and general observation that the conjugacy class of $x$ in $G$ splits into possibly several conjugacy classes in $H$ which are in one-to one correspondence with the cosets of
 $$G/ H C_G(x), $$
 where $C_G(x)$ is the centralizer of $x$ in $G$. Since $\Gamma$ is a free group, it is obvious in our case that whenever $x\neq Id$, 
 $$C_G(x)=\{ x^k\ :\ k \in  \Z \},$$
 and therefore
 $$G/ H C_G(x)=G/H.$$

 Going back to our trace computations, we have {\it formally} obtained 
 $$\det(I-\lt_s)=\mathrm{exp}\left ( -\sum_{N=1}^\infty  \frac{1}{N} \mathrm{Tr}(\lt_s^N)\right)$$
 $$=\mathrm{exp}\left( -\sum_{ \mathrm{C} \in \mathcal{P}(\Gamma(q))} \sum_{j=1}^\infty \frac{1}{j} \sum_{k=0}^\infty e^{-j(s+k)\ell(\mathcal{C})}    \right)
 =Z_{\Gamma(q)}(s),$$
 where $\mathcal{P}(\Gamma(q))$ is the set of primitive conjugacy classes in $\Gamma(q)$. To justify convergence, first observe that we have
 $$ \mathrm{Tr}(\lt_s^N)=\sum_{\ell} \sum_{\gamma_\alpha(\D_\ell)\subset \D_\ell \atop \gamma_\alpha \equiv \mathrm{Id\ mod}\ q} 
 \frac{(\gamma_\alpha'(x_\alpha))^s}{1-\gamma_\alpha'(x_\alpha)},$$
 which is roughly bounded by
 \begin{equation}
 \label{bound1}
 \vert \mathrm{Tr}(\lt_s^N) \vert \leq \sum_{B^N w=w} \frac{\left ( (B^N)'(w)\right)^{-\Re(s)}}{1-[(B^N)'(w)]^{-1}},
 \end{equation}
 and the pressure formula (\ref{pressure}) together with Bowen's result show uniform convergence of the series
 $$\sum_{N\geq 1} \frac{1}{N}  \mathrm{Tr}(\lt_s^N) $$
 on half-planes $\{\Re(s)\geq \sigma_0>\delta \}$, uniformly in $q$. Moreover, by using Proposition \ref{normest1}, we know that the {\it spectral radius} of $\lt_s$ is bounded (uniformy in $q$) by $e^{P(\Re(s))}$, therefore we do have 
 $$\det(I-\lt_s)=\mathrm{exp}\left ( -\sum_{N=1}^\infty  \frac{1}{N} \mathrm{Tr}(\lt_s^N)\right)$$
 for all $\Re(s)>\delta$. Since the infinite product formula for the Selberg's zeta function holds for all $\Re(s)>\delta$, we have the desired conclusion by analytic continuation. 
 $\square$

 \bigskip
 The above formula is critical in our analysis : a result of Patterson and Perry \cite{PatPer} says that resonances (apart from topological zeros located at negative integers) are the same (with multiplicity) as zeros of the Selberg zeta function. Therefore resonances on $\Gamma(q)\backslash \H$ are the same as non-trivial zeros of $\det(I-\lt_s)$, with multiplicities. Such a correspondence
 is also pointed out in \cite{BGS}, but comes after a rather roundabout argument based on different calculations of Laplace transforms of counting functions. 
 
\subsection{Proof of the basic upper bound}
The goal of this section is to prove the following bound.
\begin{propo}
\label{detest1}
 There exists a constant $C_\Gamma$ such that for all $q$ large and all $s \in \C$, we have the bound 
 $$\log\vert \det(I-\lt_s) \vert \leq C_\Gamma q^{3}\log (q)( 1+ \vert s \vert^2).$$
\end{propo}
\noindent The proof will follow from a careful estimate of singular values of the operators $\lt_s:H^2_q \rightarrow H^2_q$. We need first to recall some material on singular values and Weyl inequalities, our basic reference is the book of Simon \cite{Simon}. Let $\mathcal{H}_1,\mathcal{H}_2$ be two Hilbert spaces. Consider $T:\mathcal{H}_1\rightarrow \mathcal{H}_2$ a compact operator. The singular value sequence $\mu_0(T)\geq \mu_1(T)\geq \ldots \mu_k(T)$ is defined as the eigenvalue sequence of
$$\sqrt{T^* T}:\mathcal{H}_1\rightarrow \mathcal{H}_1.$$ We will need to use the following fact.

\begin{lem}
\label{vector}
Let $\mathcal{H}_1,\ldots,\mathcal{H}_m$ be $m$ Hilbert spaces with Hilbert bases denoted by 
$$(\mathrm{e}_\ell^1)_{\ell \in \N},\ldots,(\mathrm{e}_\ell^m)_{\ell \in \N}.$$
Let 
$$T=\left [T_{i,j} \right ]_{1\leq i,j\leq m}:\mathcal{H}_1\oplus \ldots \oplus \mathcal{H}_m \rightarrow \mathcal{H}_1\oplus \ldots \oplus \mathcal{H}_m$$ be a compact 
operator where each $T_{i,j}:\mathcal{H}_j \rightarrow \mathcal{H}_i$. 
Then we have for all $k\geq 0$, 
$$\mu_k(T)\leq \#\{(i,j)\  :\  T_{i,j}\neq 0 \} . \max_{i,j} \left ( \sum_{\ell \geq [k/m]} \Vert T_{i,j} \mathrm{e}_\ell^j \Vert_{\mathcal{H}_i}\right ).$$
\end{lem}
\noindent {\it Proof.}
We set $\mathcal{H}=\mathcal{H}_1\oplus \ldots \oplus \mathcal{H}_m$. We define a natural basis $(\mathbf{e}_\ell^k)_{\ell,k}$ of $\mathcal{H}$ by setting for
all $k=1,\ldots,m$
$$ \mathbf{e}_\ell^k:=(0,\ldots,0,\underbrace{\mathrm{e}_\ell^k}_{k},0,\ldots,0).$$
From the min-max principle for the eigenvalues of compact self-adjoint operators it follows that 
$$\mu_k(T)=\min_{\mathrm{codim}(V)=k} \max_{v\in V, \atop \Vert v \Vert=1} \Vert T  v \Vert_{\mathcal{H}},$$
where the min is taken along all subspaces $V\subset \mathcal{H}$ with codimension $k$.
Choosing $V=\mathrm{Span}\{ \mathbf{e}_\ell^k\ :\ \ell \geq N,\ k=1,\ldots,m\} $ we have
$$\mu_{mN}(T)\leq \max_{v\in V, \atop \Vert v \Vert=1} \Vert T  v \Vert_{\mathcal{H}}. $$
Writing 
$$v=\sum_{1\leq k \leq m,\atop \ell \geq N} \langle v,\mathbf{e}_\ell^k \rangle_{\mathcal{H}} \mathbf{e}_{\ell}^k,$$
we obtain by the triangle inequality and Cauchy-Schwarz
$$\mu_{mN}(T)\leq \sum_{1\leq k \leq m,\atop \ell \geq N} \Vert T \mathbf{e}_\ell^k\Vert_{\mathcal{H}}.$$
Let $P_j:\mathcal{H}\rightarrow \mathcal{H}$ be defined by
$$P_j(x_1,\ldots,x_m):=(0,\ldots,0,\underbrace{x_j}_{j},0,\ldots,0),$$
so that we can write
$$ T=\sum_{i,j} P_i T P_j.$$
We have obviously 
$$\Vert T \mathbf{e}_\ell^k \Vert_{\mathcal{H}}\leq \sum_{i,j} \Vert P_i T P_j \mathbf{e}_\ell^k \Vert_{\mathcal{H}},$$
hence
$$ \mu_{mN}(T)\leq \sum_{i,j} \sum_{\ell \geq N} \Vert P_i T P_j \mathbf {e}_\ell^j \Vert_{\mathcal{H}}$$
$$=\sum_{i,j} \sum_{\ell \geq N} \Vert T_{i,j}\mathrm {e}_\ell^j \Vert_{\mathcal{H}_i}\leq \#\{(i,j)\  :\  T_{i,j}\neq 0 \}
\max_{i,j} \sum_{\ell \geq N} \Vert T_{i,j}\mathrm {e}_\ell^j \Vert_{\mathcal{H}_i}.$$
The proof ends by writing
$$\mu_k(T)\leq \mu_{m[k/m]}(T)$$
and applying the above formula. $\square$

We can now move on to the proof of Proposition \ref{detest1}.  Viewing $H^2_q$ as
$$H^2_q=\bigoplus_{g \in \mathbf{G}}  H^2_g(\Omega),$$
where $\Omega=\cup_{j=1}^{2p} \D_j$,  the formula
$$\lt_s(F)(z,g):=\sum_{j\neq i} (\gamma_j'(z))^s F(\gamma_j z, \gamma_j g ),$$
shows that in the matrix representation of $\lt_s$, there are at most $2p$ non-zero operator entries per row. Using the explicit basis $(\mathbf{e}_k^\ell)$ for $H^2(\Omega)$,
it is enough to estimate 
$$\Vert (\gamma_j')^s \mathbf{e}_k^\ell\circ \gamma_j \Vert_{H^2(\D_i)}, $$
where $\gamma_j(\D_i)\subset \D_\ell$. Using the fact that 
 $$\mathrm{dist}(\gamma_j(\D_i),\partial \D_\ell )\geq \epsilon_0,$$
 we obtain the bound (for some adequate constants $C>0$ and $0< \rho_0<1$ )
 $$ \Vert (\gamma_j')^s \mathbf{e}_k^\ell\circ \gamma_j \Vert_{H^2(\D_i)}\leq Ce^{C\vert s \vert} \rho_0^k.$$
 Applying Lemma \ref{vector}, we have reached 
 $$\mu_k(\lt_s)\leq 2p \vert G \vert C e^{C\vert s\vert} \sum_{j\geq [ k/\vert G \vert ]} \rho_0^j \leq \widetilde{C} \vert G \vert e^{C \vert s \vert } \rho_0^{\frac{k}{\vert G \vert}}.$$ 
 We can now use Weyl inequalities (see \cite{Simon} Theorem 1.15) to write
 $$\log \vert \det(I-\lt_s) \vert  \leq \sum_{k=0}^\infty \log(1+\mu_k(\lt_s) )$$
 $$\leq N \log(1+\widetilde{C}\vert G \vert e^{C\vert s \vert} )+\widetilde{C} \vert G \vert  e^{C\vert s \vert } \sum_{k>N}  \rho_0^{\frac{k}{\vert G \vert}}.$$
 Setting 
 $$N=\frac{C\vert G\vert \vert s \vert } {\vert \log \rho_0 \vert},$$
 we end up with 
 $$ \log \vert \det(I-\lt_s) \vert  \leq C' (  \vert G \vert \vert s \vert^2+\vert s \vert \vert G \vert \log \vert G \vert  +1),$$
 for some suitable constant $C'$. The proof is done since $\vert G \vert \asymp q^3$. $\square$ 
 \begin{cor}
 \label{countingest1}
 There exists a constant $C>0$ such that for all $q$ large enough, we have
 $$\mathcal{N}_q(r)\leq Cq^3\log(q)(1+r^2).$$
 \end{cor}
 \noindent {\it Proof}. This estimate follows straightforwardly from Jensen's formula (see the end of $\S 4$, Proposition \ref{Jensen} for details ) but a lower bound is required. Indeed, using the bound (\ref{bound1}), we observe that
 $$\log\vert Z_{\Gamma(q)}(1)\vert \geq -\sum_{N=1}^\infty \frac{1}{N} \sum_{B^N w=w} \frac{\left ( (B^N)'(w)\right)^{-1}}{1-[(B^N)'(w)]^{-1}}> -\infty$$
which is a lower bound independent of $q$. Applying the classical Jensen's identity on the disc $D(1, \widetilde{r})$, where $\widetilde{r}$ is carefully chosen in function of $r$, we end up with the above bound. $\square$
 
\section{Nielsen volume and trace formula, proof of the first theorem}
In this section, we use the global upper bound proved in the previous section to produce a lower bound, thanks to the leading singularity of the trace formula.
Before we state the trace formula, we need to point out a fact. The Nielsen volume of a geometrically finite surface $X$ is defined as the hyperbolic area of the {\it Nielsen region }
$N$, the geodesically convex hull of closed geodesics on the surface. In the convex co-compact case, the Nielsen region is a compact surface with geodesic boundary. 
From the Gauss-Bonnet formula, we know (see for example \cite{Borthwick} Theorem 2.15) that $\mathrm{Vol}(N)=-2\pi \chi(N)$, where $\chi(N)=\chi(X)$ is the Euler-Poincar\'e characteristic.
Going back to our surfaces $X_q=\Gamma(q)\backslash \H$, there is a natural covering
$$\Gamma(q)\backslash \H \rightarrow \Gamma \backslash \H,$$
with degree $[\Gamma(q):\Gamma]=\vert G \vert$. It is a standard fact of algebraic topology that $\chi(X_q)=\vert G \vert \chi(X)$ which translates into the formula
$$\mathrm{Vol}(N_q)=\vert G \vert \mathrm{Vol}(N)\asymp q^3.$$

\noindent 
The Wave-trace formula stated below is due to Guillop\'e and Zworski
\cite{GuiZwor2}. We denote by ${\mathcal P}_q$ the set of {\it primitive closed geodesics} on the surface $X=\Gamma(q)\backslash \H$, 
and if $\gamma\in {\mathcal P}_q$, $l(\gamma)$
is the length. In the following, $N_q$ still denotes the Nielsen region. Let $\varphi\in C_0^\infty( (0,+\infty))$ i.e. a smooth function, compactly supported in $\R_+^*$. We have the identity:
\begin{equation}
\label{selberg}
\sum_{s\in {\mathcal R}_q}\widehat{\varphi}(i(s-\half))=
-\frac{{\rm Vol}(N_q)}{4\pi}\int_0^{+\infty} \frac{\cosh(t/2)}{\sinh^2(t/2)}\varphi(t)dt.
\end{equation}
$$+\sum_{\gamma\in {\mathcal P}_q}\sum_{k\geq 1}\frac{l(\gamma)}{2\sinh(kl(\gamma)/2)}\varphi(kl(\gamma)),$$
where $\widehat{\varphi}$ is the usual Fourier transform
$$\widehat{\varphi}(\xi)=\int_{\R}\varphi(x)e^{-ix\xi}dx.$$
To prove Theorem \ref{lowerb}, we choose a test function $\varphi_0 \in C_0^\infty(0,1)$ such that $\varphi_0 \geq 0$ and $\int \varphi_0=1$. For all $T>0$ we set 
$$\varphi_T(x):=T\varphi_0(T x), $$
where $T$ will be a large parameter. Since the length spectrum of $X_q$ is a subset of the length spectrum of $X=\Gamma \backslash \H$ (without multiplicities), we can definitely find a uniform $\epsilon_0>0$ such that for all
$\gamma \in \mathcal{P}_q$, $l(\gamma)\geq \epsilon_0$. In the sequel, we take $T$ large enough so that $T^{-1}<\epsilon_0$. The above trace formula gives
$$\sum_{s\in {\mathcal R}_q}\widehat{\varphi_T}(i(s-\half))=
-\frac{{\rm Vol}(N_q)}{4\pi}\int_0^{+\infty} \frac{\cosh(t/2)}{\sinh^2(t/2)}\varphi_T(t)dt.
 $$
In view of the preceding remarks, this yields for all $q$ large and $T$ as above
$$\left \vert \sum_{s\in {\mathcal R}_q}\widehat{\varphi_T}(i(s-\half))  \right \vert  \geq C q^3 T^2,$$
where $C>0$ is a uniform constant. Using the fact that 
$$\widehat{\varphi_T}(\xi)=\widehat{\varphi_0}(\xi/T)$$
 and repeated integrations by parts, we have the estimate
$$\vert \widehat{\varphi_T}(z) \vert \leq C_N \frac{e^{\Im(z)/T}}{(1+\vert z/T \vert )^N},$$
for all $N\geq 0$. We write 
$$\left \vert \sum_{s\in {\mathcal R}_q}\widehat{\varphi_T}(i(s-\half))  \right \vert\leq
\sum_{\vert s \vert \leq R}\vert \widehat{\varphi_T}(i(s-\half)) \vert  +\sum_{\vert s \vert> R}\vert \widehat{\varphi_T}(i(s-\half)) \vert$$
$$\leq A \mathcal{N}_q(R)+C_N \sum_{\vert s \vert >R} \frac{1}{\left (1+\frac{\vert s-1/2 \vert}{T} \right)^N},$$
where $A>0$ does not depend on $q,T,R$. Since we have obviously 
$$\vert s-1/2\vert \geq \vert s \vert -1/2,$$
and $T$ is taken large, we can write
$$\sum_{\vert s \vert >R} \frac{1}{\left (1+\frac{\vert s-1/2 \vert}{T} \right)^N}\leq \int_R^\infty \frac{d \mathcal{N}_q(t)}{(t/T)^N}.$$
A Stieltjes integration by parts yields
$$\int_R^\infty \frac{d \mathcal{N}_q(t)}{(t/T)^N}\leq \frac{\mathcal{N}_q(R)}{(R/T)^N}+N \int_{R/T}^\infty\frac{\mathcal{N}_q(xT)dx}{x^{N+1}}$$
$$\leq C_1 (q^3\log q) R^2 \left( \frac{R}{T}\right)^{-N}+\widetilde{C_N} (q^3\log q)T^2 \left(\frac{R}{T} \right)^{2-N},$$
where we have used the upper bound from Corollary \ref{countingest1}.

\bigskip
Setting $$R=T(\log q)^\epsilon,$$ with $\epsilon>0$ and $T$ large (but fixed), we have obtained
$$Cq^3T^2\leq A\mathcal{N}_q(R)+ C_N' q^3T^2 (\log q)^{1+(2-N)\epsilon},$$
where $C'_N$ is a (possibly large) constant depending only on $N$. Taking $N$ so large that $1+(2-N)\epsilon <0$, we get that for all $q$ large enough,
$$C'_N(\log q)^{1+(2-N)\epsilon}\leq \frac{C}{2},$$
which yields
$$A^{-1}\frac{C}{2}T^2 q^3 \leq \mathcal{N}_q( T (\log q)^\epsilon),$$
and the proof is done. $\square$

\section{The refined upper bound and proof of the second theorem} 

\subsection{Refined function space $H^2_q(h)$ } 
Let $0<h$ and
set $$\Lambda(h):=\Lambda(\Gamma)+(-h,+h),$$
then for all $h$ small enough, $\Lambda(h)$ is a bounded subset of $\R$ whose connected components have length at most $Ch$ where $C>0$ is independent of $h$, see \cite{Borthwick} Lemma 15.12. Let $\{I_\ell(h),\ \ell=1,\ldots,N(h)\}$ denote these connected components. The existence of a finite Patterson-Sullivan measure 
$\mu$ supported on $\Lambda(\Gamma)$  shows that (see \cite{Borthwick} P. 312, first displayed equation)
the number $N(h)$ of connected components is $O\left(h^{-\delta} \right)$. Given $1\leq \ell \leq N(h)$, let $\D_\ell(h)$ be the unique euclidean open disc
in $\C$ orthogonal to $\R$ such that $$\D_\ell(h)\cap \R=I_\ell(h).$$ 
We consider
$$H^2_q(h):= \bigoplus_{\ell=1}^{N(h)}H^2(\D_\ell(h)\times \mathbf{G}).$$ 
Set $$\Omega(h):=\cup_{\ell=1}^{N(h)}\D_\ell(h),$$ then the norm on $H^2_q(h)$ is given by 
$$ \Vert F \Vert_{q,h}^2:=\sum_{g \in \G} \int_{\Omega(h)} \vert F(z,g)\vert^2dm(z).$$
The parameter $h$ will play the role of a scale parameter whose size will be adjusted according to the spectral parameter $s$.
An important fact in the sequel is the following estimate taken from \cite{NaudInvent}.
\begin{lem}
 \label{separation}
 There exists $n_0$ such that for all $n\geq n_0$, for all $\alpha\in \mathscr{W}_n^j$ and all 
 $\ell \in \mathscr{E}_j(h)$, there exists an index $\ell'$ such that $\gamma_\alpha(\D_\ell(h))\subset \D_{\ell'}(h)$ with
 $$ \mathrm{dist}(\gamma_\alpha(\D_\ell(h)),\partial \D_{\ell'}(h) )\geq \half h.$$
\end{lem}
The above Lemma guarantees that the transfer operator considered previously is well defined for all $n$ large enough (independently of $h,q$):
$$\lt_s^n:H^2_q(h)\rightarrow H^2_q(h).$$
The basic norm estimate is the following.
\begin{propo}
 \label{normest2}
 Set $\sigma=\Re(s)$, where $s$ is the spectral parameter.
 There exist a constant $C_\sigma>0$, independent of $q,h$ such that for all $n$ large and 
 $$\Vert \lt_s^n \Vert_{H^2_q(h)}\leq C_\sigma e^{C_\sigma h \vert \Im(s)\vert} h^{-\delta} e^{nP(\sigma)} .$$
\end{propo}
The proof is postponed to the appendix. This estimate essentially shows that the spectral radius of $\lt_s^n$ can be bounded uniformly on all spaces $H^2_q(h)$ in term of the topological pressure. Moreover, by following verbatim the trace computations of $\S 2$, we see that the determinants (and hence the full spectrum of $\lt_s^n$) 
$$\det(I-z\lt_s^n)$$
{\it do not} depend on the scale parameter $h$ (traces depend only on periodic points of the Bowen-Series map). To count resonances on $\mathcal{R}_q$, we will
use the following family of {\it Hilbert-Schmidt}  determinants
$$\zeta_{n}(s):={\det}_2 ( I-\lt^n_s).$$
Remark that if $s\in \mathcal{R}_q$, then by Proposition \ref{detformula}, the operator 
$$\lt_s:H^2(q) \rightarrow H^2(q)$$ must have $1$ as an eigenvalue, and so does
$\lt_s^n$, but clearly $H^2(q)\hookrightarrow H^2_q(h)$ for all $h$ small, therefore $\zeta_{n}(s)=0$. On the other hand, $\zeta_{n}(s)$ might have some extra zeros which are not resonances, but that's a minor issue since we are interested on upper bounds on the density.
\subsection{Two observations and a consequence}
In this subsection, we prove to lemmas which are to be used in the proof of the main estimate in a crucial way. However, since they only play a role at the very end of the proof, the reader can skip them at first glance.
\begin{lem}
\label{girth}
There exists $\epsilon_1(\Gamma)>0$ such that for all $j=1,\ldots,2p$, all $\alpha,\beta \in {\mathscr W}_n^j$ with 
$\gamma_\alpha \equiv \gamma_\beta\ \mathrm{mod}\ q$, we have 
$$n<\epsilon_1 \log q \Rightarrow \alpha=\beta. $$
\end{lem}
\noindent {\it Proof.} This is a slight variation on the "girth lower bound" proved in \cite{Gamburd} for Cayley graphs of $SL_2(\F_q)$ with respect to arbitrary generators of $\Gamma$. Let $\Vert .\Vert$ be the usual enclidean norm on $\R^2$ and if $M$ is a $2\times 2$ real matrix, set
$$\Vert M \Vert=\sup_{\Vert X\Vert \leq 1} \Vert MX \Vert,$$
which is an algebra norm. Given $j\in \{1,\ldots,2p \}$, assume that we have two words $\alpha,\beta  \in {\mathscr W}_n^j$ with 
$$\gamma_\alpha \equiv \gamma_\beta\ \mathrm{mod}\ q\ \mathrm{and}\ \gamma_\alpha \neq \gamma_\beta.$$
Consider the matrix $\gamma_\alpha \gamma_{\beta}^{-1}$, then we do have
$$\gamma_\alpha \gamma_\beta^{-1}\equiv \left (\begin{tabular}{cc}1&0\\0&1 \end{tabular}\right)\  \mathrm{mod}\ q,$$
but 
$$\gamma_\alpha \gamma_\beta^{-1}\neq \left (\begin{tabular}{cc}1&0\\0&1 \end{tabular}\right).$$
Therefore one of the two off-diagonal entries of $\gamma_\alpha\gamma_\beta^{-1}$ is a non-zero multiple of $q$, which forces
$$\Vert \gamma_\alpha \gamma_\beta^{-1} \Vert \geq q.$$
Since $\Vert. \Vert$ is an algebra norm, we have 
$$q\leq \left (\max_{j=1,\ldots,2p}  \Vert \Gamma_j \Vert \right)^{2n} $$
and the proof is done with
$$\epsilon_1=\left (\max_{j=1,\ldots,2p}  \Vert \Gamma_j \Vert \right)^{-2}.\ \square$$
We will also need to use the following fact.
\begin{lem}
\label{separate}
There exist constants $\overline{C}>0$ and $0<\overline{\theta}<1$ such that for all $j=1,\ldots,2p$, for all $z\in \D_j$ and all words $\alpha\neq \beta \in \mathscr{W}_n^j$,
$$\vert \gamma_\alpha(z)-\gamma_\beta(z) \vert \geq \overline{C}\overline{\theta}^{r(\alpha,\beta)} ,$$
where $r(\alpha,\beta)=\max\{ 0\leq i \leq n\ :\ \forall k\leq i,\  \alpha_k=\beta_k \}$.
\end{lem}
\noindent {\it Proof}. Let $\alpha\neq \beta \in \mathscr{W}_n^j$ and pick $z\in \D_j$. Since $\alpha\neq \beta$ we have $r(\alpha,\beta)\leq n-1$. Let us write
$$\vert \gamma_\alpha(z)-\gamma_\beta(z)\vert=\vert \widetilde{\gamma}(w_1)-  \widetilde{\gamma}(w_2)  \vert $$
where 
$$\widetilde{\gamma}(w)=\gamma_{\alpha_1}\circ \gamma_{\alpha_2}\circ \ldots \gamma_{\alpha_r}(w)=
\gamma_{\beta_1}\circ \gamma_{\beta_2}\circ \ldots \gamma_{\beta_r}(w),$$
and $w_1=\gamma_{\alpha_{r+1}}\circ \ldots\circ \gamma_{\alpha_n}(z)$, $w_2=\gamma_{\beta_{r+1}}\circ \ldots\circ \gamma_{\beta_n}(z)$.
Since $\widetilde{\gamma}$ is a M\"obius transform, we can use the standard formula 
$$\vert \widetilde{\gamma}(w_1)-\widetilde{\gamma}(w_2)\vert^2=\vert\widetilde{\gamma}'(w_1) \vert \vert \widetilde{\gamma}'(w_2)\vert \vert w_1-w_2 \vert^2.$$
Recall that for all $k=1,\ldots,2p$, for all $i\neq k$ 
$$\gamma_k(\D_i)\subset \D_{p+k},$$
where $p+k$ is understood mod $2p$. Therefore $\alpha_{r+1}\neq \beta_{r+1}$ implies that $w_1$ and $w_2$ belong to two different discs
$$ w_1 \in \D_{p+\alpha_{r+1}}\neq \D_{p+\beta_{r+1}}  \ni w_2. $$
Therefore we have 
$$\vert w_1-w_2\vert\geq \min_{k\neq \ell} \mathrm{dist}( \D_k,\D_\ell) >0.$$
Using the lower bound for the derivatives from (\ref{hyp1}), we end up with
$$\vert \gamma_\alpha(z)-\gamma_\beta(z)\vert\geq  \min_{k\neq \ell} \mathrm{dist}( \D_k,\D_\ell) C^{-1} \overline{\theta}^{r(\alpha,\beta)},$$
and the proof is done. $\square$

\bigskip
Both of these estimates are of independent interest but we will actually combine them as follows.
\begin{cor}
\label{Erhenfest} 
Let $C>0$ be a constant. There exists $\epsilon_0>0$ depending only on $C,\Gamma$ such that for all $j=1,\ldots,2p$, for all $z\in \D_j$ and all $\alpha,\beta \in \mathscr{W}_n^j$
with $n\leq \epsilon_0 (\log q +\log h^{-1})$, we have 
$$\gamma_\alpha \equiv \gamma_\beta\ \mathrm{mod}\ q\ \mathrm{and}\ \vert \gamma_\alpha(z)-\gamma_\beta(z) \vert \leq C h\ \Rightarrow \alpha=\beta.$$
\end{cor}
\noindent {\it Proof}.
Let $n\leq \epsilon_0 (\log q + \log h^{-1})$ where $\epsilon_0$ will be adjusted later on. Assume that we have two words $\alpha\neq \beta \in \mathscr{W}_n^j$ such that
$$ \gamma_\alpha \equiv \gamma_\beta\ \mathrm{mod}\ q\ \mathrm{and}\ \vert \gamma_\alpha(z)-\gamma_\beta(z) \vert \leq C h.$$
By Lemma \ref{separate}, we get that 
$$\overline{C} \overline{\theta}^{n-1}\leq \vert \gamma_\alpha(z)-\gamma_\beta(z) \vert \leq C h, $$
which shows that 
$$\log( h^{-1})\leq n\log (\overline{\theta}^{-1} )+\widetilde{C},$$
where $\widetilde{C}$ is another constant (depending on the previous ones).  Assuming 
$$\epsilon_0 \log(\overline{\theta}^{-1} )<1,$$ we get
$$n\leq \frac{\epsilon_0}{1-\epsilon_0 \log(\overline{\theta}^{-1} )} \log q + C'.$$
It is now clear that if $\epsilon_0$ is taken small enough, we have for large $q$
$$n<\epsilon_1 \log q, $$
hence contradicting Lemma \ref{girth}. $\square$
\subsection{Hilbert-Schmidt norms and pointwise estimate} 
The main theorem will follow, after a suitable application of Jensen's formula from the next statement which is the main goal of this section. We recall that we will work with
the modified zeta function
$$\zeta_{(n)}(s):=\mathrm{det}_{2}(I-\lt_s^n),$$
where $n=n(q,T)$ will be adjusted later on. 

\begin{propo}
\label{firstest}
Fix $\delta>\sigma>\delta/2$. Then there exist constants $\epsilon_0>0$, $C_\sigma>0$ and $\eta_j(\sigma)>0$, $j=1,2$ such that for all 
$0\leq \vert \Im(s) \vert \leq T$ (with $T\geq 1$) and $\sigma\leq \Re(s) \leq \delta$, we have for all $q$ large, 
$$\log \vert \zeta_{n(T,q)}(s) \vert \leq C_\sigma T^{\delta-\eta_1(\sigma)} q^{3-\eta_2(\sigma)},$$
where $n(T,q)=[\epsilon_0( \log q +\log h^{-1})]$, $h=T^{-1}$.
\end{propo}
It is necessary to recall at this point a few facts on regularized Hilbert-Schmidt determinants, our reference is \cite{Gohberg}.
Let $\mathcal{H}$  be an abstract separable Hilbert space, and 
$T:\mathcal{H}\rightarrow \mathcal{H}$ a compact operator. The operator $T$ is called Hilbert-Schmidt if 
$$\sum_{k=0}^\infty \mu_k(T)^2<\infty,$$
and its Hilbert-Schmidt norm is 
$$\Vert T \Vert_{HS}^2:=\mathrm{Tr}(T^* T)=\sum_{k=0}^\infty \mu_k(T)^2. $$
The regularized determinant ${\det}_2(I+T)$ is defined for all Hilbert-Schmidt operators by 
$${\det}_2(I+T):=\det(I+ [(I+T)\mathrm{exp}(-T)-I]),$$
where $(I+T)\mathrm{exp}(-T)-I$ is a trace class operator. If $T$ is itself a trace class operator, then we have actually
$${\det}_2(I+T)=\det(I+T) e^{-\mathrm{Tr}(T)}.$$
A key tool for our purpose is the following inequality \footnote{It is also possible to work with the usual Fredholm determinants: one has to consider instead
$\zeta_{n}(s):=\det(I-\lt_s^{2n})$ and use the inequality
$$\log\vert \det(I-T^{2n}) \vert \leq \Vert T^{2n}\Vert_{Tr}\leq \Vert T^n \Vert_{HS}^2.$$}
 (see \cite{Gohberg}, chapter 4, Theorem 7.4) :
\begin{equation}
\label{HSest}
\vert {\det}_2(I+T) \vert \leq e^{\half \Vert T \Vert^2_{HS}}.
\end{equation}

We can now give a proof of Proposition \ref{detest1}. First we will use the notation for all $j=1,\ldots,2p$
$$\Omega_j(h)=\Omega(h)\cap \D_j.$$ 
Given $\ell \in \{1,\ldots,N(h)\}$, let $(\mathbf{e}_k^\ell)_{k\in \N}$ be a Hilbert basis of $H^2(\D_\ell(h))$.

According to inequality (\ref{HSest}), we need to estimate the Hilbert-Schmidt norm
$$\Vert \lt_s^n \Vert_{HS}^2=\mathrm{Tr}((\lt_s^n)^* \lt_s^n)$$
$$=\sum_{g \in \G} \sum_{k,\ell} \sum_{w\in \G} \int_{\Omega(h)} \vert \lt_s^n( \mathbf{e}_k^\ell \otimes \mathscr{D}_g)(z,w)\vert^2 dm(z).$$
In addition, we have
$$\int_{\Omega(h)} \vert \lt_s^n( \mathbf{e}_k^\ell \otimes \mathscr{D}_g)(z,w)\vert^2 dm(z)$$
$$=\sum_{j=1}^{2p} \sum_{\alpha,\beta \in \mathscr{W}_n^j} 
\mathscr{D}_g(\gamma_\alpha w) \mathscr{D}_g(\gamma_\beta w) \int_{\Omega_j(h)} (\gamma_\alpha')^s \overline{(\gamma_\beta')^s}  \mathbf{e}_k^\ell\circ \gamma_\alpha \overline{\mathbf{e}_k^\ell\circ \gamma_\beta} dm.$$
Noticing that 
$$\sum_{g \in \G} \sum_{w \in \G} \mathscr{D}_g(\gamma_\alpha w) \mathscr{D}_g(\gamma_\beta w)=\left \{ \vert \G \vert\ \mathrm{if}\ \gamma_\alpha \equiv 
\gamma_\beta\ [q] \atop 0\ \mathrm{otherwise} \right.$$
We have obtained 
$$\Vert \lt_s^n \Vert_{HS}^2=\vert \G \vert \sum_{k,\ell} \sum_{j=1}^{2p} \sum_{\alpha,\beta \in \mathscr{W}_n^j \atop \gamma_\alpha \equiv \gamma_\beta\ [q]} 
 \int_{\Omega_j(h)} (\gamma_\alpha')^s \overline{(\gamma_\beta')^s}  \mathbf{e}_k^\ell\circ \gamma_\alpha 
 \overline{\mathbf{e}_k^\ell\circ \gamma_\beta} dm.$$
 Since 
 $$\sum_{k,\ell} \mathbf{e}_k^\ell(z_1)
 \overline{\mathbf{e}_k^\ell(z_2)}  $$
 converges uniformly on compact subsets of $\Omega(h)\times \Omega(h)$ to the Bergman kernel $B_{\Omega(h)}(z_1,z_2)$, we can exchange summations to write
 \begin{equation}
 \label{Magicidentity}
  \Vert \lt_s^n \Vert_{HS}^2= 
  \vert \G \vert \sum_{j=1}^{2p} \sum_{\alpha,\beta \in \mathscr{W}_n^j \atop \gamma_\alpha \equiv \gamma_\beta\ [q]} 
 \int_{\Omega_j(h)} (\gamma_\alpha'(z))^s \overline{(\gamma_\beta'(z))^s}  B_{\Omega(h)}(\gamma_\alpha z,\gamma_\beta z) dm(z). 
 \end{equation}
 We stress that since $\Omega(h)$ is disconnected, we have 
 $B_{\Omega(h)}(z,w)=0$ if $z$ and $w$ do not belong to the same connected component.
 We assume from now on that $h=T^{-1}$ with $\vert \Im(s) \vert \leq T$ and $\halfdelta <\sigma\leq \Re(s) \leq \delta$. We will choose $n:=n(q,T)$ of the form
 $$n(q,T)=[\epsilon_0 ( \log q +\log T)].$$
 We recall that each disc $\D_\ell(h)$ has by construction diameter at most $Ch$, and we choose $\epsilon_0(C)$ so that the conclusion of Corollary \ref{Erhenfest} is true.
 Therefore, in the above sum, there is no off-diagonal contribution. Indeed, according to Corollary \ref{Erhenfest} there are no words with $\alpha\neq \beta$ such that $$ \gamma_\alpha \equiv \gamma_\beta\ \mathrm{mod}\ q\ \mathrm{and}\ \vert \gamma_\alpha(z)-\gamma_\beta(z) \vert \leq C h$$
provided that $n\leq \epsilon_0(C)(\log q+ \log(h^{-1})$.  As a consequence we have actually
\begin{equation}
 \label{Magicidentityprim}
  \Vert \lt_s^{n(q,T)} \Vert_{HS}^2= 
  \vert \G \vert \sum_{j=1}^{2p} \sum_{\alpha \in  \mathscr{W}_n^j } 
 \int_{\Omega_j(h)} \vert (\gamma_\alpha'(z))^s \vert^2  B_{\Omega(h)}(\gamma_\alpha z,\gamma_\alpha z) dm(z). 
 \end{equation}

 Using the fact \footnote{The size of each disc compensates exactly for the exponential growth of $(\gamma'_\alpha)^s$ as $\Im(s)$ becomes large, see \cite{NaudInvent}, after Lemma 3.4, P.737.}
  that each disc $\D_\ell(h)$ has radius at most $Ch=CT^{-1}$, and because of the uniform distortion estimate (\ref{bdist}), we have for all $\alpha \in \mathscr{W}_n^j$ and all $z\in \Omega_j(h)$,
 $$\vert (\gamma'_\alpha(z))^s \vert \leq C' \sup_{I_j} (\gamma'_\alpha)^{\Re(s)}.$$
 On the other hand, using Lemma \ref{separation} and the explicit formula for the Bergman Kernel, we see that
 $$\vert B_{\Omega(h)}( \gamma_\alpha z, \gamma_\alpha z) \vert \leq C''h^{-2},$$
uniformly in $n$. Since we have
$$m(\Omega(h))=O(h^{2-\delta}),$$
we obtain by the inequality (\ref{pressure2}) that
$$ \Vert \lt_s^n \Vert_{HS}^2\leq C_\sigma \vert G \vert h^{-\delta} e^{n(q,T)P(2\sigma)}\leq C'_\sigma q^3 T^\delta e^{n(q,T)P(2\sigma)}.$$
Now recall that because of Bowen's formula, $\sigma>\delta/2$ implies that $P(2\sigma)<0$ and the proof is done with 
$$\eta_1(\sigma)=\eta_2(\sigma)=-\epsilon_0 P(2\sigma)$$
since $n(q,T)\geq \epsilon_0(\log q +\log T)-1$. $\square$
\subsection{Applying Jensen's formula} 
Using Proposition \ref{detest1}, we can prove Theorem \ref{upperb}. We will apply the following version of Jensen's formula which can be derived straightforwardly from
the classical textbooks, for example \cite{Tit}, p.125.
\begin{propo}
\label{Jensen}
   Let $f$ be a holomorphic function on the open disc $D(w,R)$, and assume that $f(w)\neq 0$. let $N_f(r)$ denote the number of zeros of $f$ in the closed disc
   $\overline{D}(w,r)$. For all $\widetilde{r}<r<R$, we have 
   $$N_f(\widetilde{r})\leq \frac{1}{\log(r/\widetilde{r})} \left (  \frac{1}{2\pi} \int_0^{2\pi} \log \vert f(w+re^{i\theta})\vert d\theta-\log\vert f(w)\vert   \right).$$
\end{propo}
The goal is to apply the above formula to $\zeta_n(s)$ where $n$ is taken according to Proposition \ref{detest1}. We need a lower bound. Assume that $\Re(s)\geq 1$,
we get
$$\zeta_n(s)=\det(I-\lt_s^n) e^{\mathrm{Tr}(\lt_s^n)}$$
$$=\mathrm{exp}\left(   -\sum_{N=1}^\infty \frac{1}{N} \mathrm{Tr}(\lt_s^{nN}) +\mathrm{Tr}(\lt_s^n) \right).$$
Using the bound (\ref{bound1}) we do have (recall that $B$ is the Bowen-Series map on the boundary).
$$\vert \mathrm{Tr}(\lt_s^{nN}) \vert \leq \sum_{B^{nN} w=w} \frac{\left ( (B^{nN})'(w)\right)^{-1}}{1-[(B^{nN})'(w)]^{-1}}.$$
On the other hand, formula (\ref{pressure}) for the topological pressure gives us (for all $\epsilon>0$)
$$\log \vert\zeta_n(s)\vert \geq -C_\epsilon\sum_{N=1}^\infty \frac{1}{N} e^{nN(P(1)+\epsilon)} -C_\epsilon e^{n(P(1)+\epsilon)}.$$
Since $P(1)<0$, this last lower bound shows clearly that one can find $\kappa>0$ independent of $s,n$ such that for all $\Re(s)\geq 1$, we have 
$$\log \vert\zeta_n(s)\vert \geq-\kappa. $$
Going back to the proof of Theorem \ref{upperb}, fix now $\halfdelta< \sigma_2<\sigma_1<\sigma_0<\delta$.  Let $\mathscr{R}(\sigma_0,T)$ denote the (closed) rectangle
$$\mathscr{R}(\sigma_0,T):=\{\sigma_0\leq \Re(s)\leq \delta\ \mathrm{and}\ \vert \Im(s)-T\vert \leq 1 \}.$$
For $r\geq 1$, set $M(r)=r-\sqrt{r^2-1} \asymp \frac{1}{r}$, and choose $r$ large enough so that $\sigma_0-\sigma_1=M(r)$.
Set $w=\sigma_1+r+iT$. Clearly  if $\sigma_1-\sigma_2$ is small enough, we do have $\Re(w)\geq 1$. One can also check that we have
$$\mathscr{R}(\sigma_0,T)\subset \overline{D}(w,r)\subset D(w,r+\sigma_1-\sigma_2).$$
Applying the above formula to $\zeta_n(s)$ with $n=n(q,T+r+\sigma_1-\sigma_2)$ on the disc 
$$D(w,r+\sigma_1-\sigma_2),$$ 
we get 
$$\mathcal{M}_q(\sigma_0,T)\leq N_{\zeta_n}(r)$$
$$\leq \frac{1}{\log(\sigma_1-\sigma_2)}\left( \frac{1}{2\pi}\int_0^{2\pi} \log \vert \zeta_n(w+(r+\sigma_1-\sigma_2)e^{i\theta})\vert d\theta   +\kappa \right).$$
Using Proposition \ref{detest1} we get 
$$\mathcal{M}_q(\sigma_0,T)\leq \frac{r+\sigma_1-\sigma_2}{\log(\sigma_1-\sigma_2)}C_{\sigma_1} (T+r+\sigma_2-\sigma_1)^{\delta-\eta_1(\sigma_2)} q^{3-\eta_2(\sigma_2)}$$
$$+\frac{\kappa}{\log(\sigma_1-\sigma_2)}.$$
Since $r,\sigma_0,\sigma_1,\sigma_2$ are fixed, we clearly get the desired conclusion, up to a change of constants.  $\square$
\subsection{Final remarks} 
Clearly the proof we have used (based on the separation Lemma \ref{separate}) not only simplifies part of the arguments in \cite{NaudInvent} but slightly strengthens the result. Moreover, we believe the technique can be carried over
to higher dimensional settings, at least for Schottky groups. On the other hand, it is clear that this trick fails for situations where the limit set is not disconnected, for example for quasi-fuchsian groups. This is where the more sophisticated
arguments (based on off-diagonal cancellations) used in \cite{NaudInvent} can be usefull, this will be pursued elsewhere. We also point out that our method should work without great modifications
to deal with subgroups of arithmetic co-compact fuchsian groups, and also for Schottky subgroups of $SL_2(\Z[i])$. The only thing required is an appropriate "logarithmic girth" estimate with respect to the congruence parameter.  
\section{Appendix : Basic norm estimates on $H^2_q(h)$}
In this section, we prove Proposition \ref{normest2}, which gives a crude bound for the operator norm of $\lt_s^N$ on spaces $H^2_q(h)$. Let $F \in H^2_q(h)$.
We first write
$$\Vert \lt_s^N(F)\Vert^2_{H^2_q(h)}=\sum_{g \in \G} \int_{\Omega(h)}\vert \lt_s^N(F)(z,g) \vert^2 dm(z)$$
$$=\sum_{g\in \G} \sum_j \sum_{\alpha,\beta \in \mathscr{W}_N^j} \int_{\Omega_j(h)} (\gamma_\alpha')^s  \overline{(\gamma_\beta')^s} 
F(\gamma_\alpha z,\gamma_\alpha g) \overline{F}(\gamma_\beta z,\gamma_\beta g)dm(z). $$
We recall that because each disc $\D_\ell(h)$ has size at most $Ch$ and thanks to the bounded distortion property, we do have
$$\sup_{z \in \D_\ell(h)} \vert (\gamma'_\alpha(z))^s \vert \leq e^{C \vert \Im(s) \vert h} \sup_{z \in \Omega_j(h)} \vert \gamma_\alpha'(z) \vert ^{\Re(s)}.$$
Therefore we have
$$\Vert \lt_s^N(F)\Vert^2_{H^2_q(h)}\leq e^{C \vert \Im(s) \vert h}  \sum_{g \in \G} \sum_j \sum_{\alpha,\beta} 
\sup_{z \in \Omega_j(h)} \vert \gamma_\alpha'(z) \vert ^{\Re(s)}\sup_{z \in \Omega_j(h)} \vert \gamma_\beta'(z) \vert ^{\Re(s)}$$ 
$$\times \int_{\Omega_j(h)} \vert F(\gamma_\alpha z, \gamma_\alpha g )\vert \vert F(\gamma_\beta z, \gamma_\beta g )\vert dm(z).$$
By the reproducing property of Bergman kernels, we have for all $z\in \Omega(h)$,
$$F(z,g)=\int_{\Omega(h)}F(w) B_{\Omega(h)}(z,w) dm(w),$$
 which allows us to write (thanks to Cauchy-Schwarz inequality and Lemma \ref{separation}) 
 $$\sup_{z \in \Omega_j(h)} \vert F(\gamma_\alpha z, g)\vert \leq C h^{-2}\sqrt{ m(\Omega(h))}  \left( \int_{\Omega(h)} \vert F(w,g)\vert^2 dm(w) \right)^{1/2}.$$
Therefore we have 
$$\int_{\Omega_j(h)} \vert F(\gamma_\alpha z, \gamma_\alpha g )\vert \vert F(\gamma_\beta z, \gamma_\beta g )\vert dm(z)
\leq C h^{-4} m(\Omega(h))^2 $$
$$\times \left( \int_{\Omega(h)} \vert F(w,\gamma_\alpha g)\vert^2 dm(w) \right)^{1/2} \left( \int_{\Omega(h)} \vert F(w,\gamma_\beta g)\vert^2 dm(w) \right)^{1/2}.$$ 
Since $m(\Omega(h))=O(h^{2-\delta})$, we have obtained
$$\Vert \lt_s^N(F)\Vert^2_{H^2_q(h)}\leq C h^{-2\delta}e^{C \vert \Im(s) \vert h}  \sum_{g \in \G} \sum_{j,\alpha,\beta} 
\sup \vert \gamma_\alpha'\vert ^{\Re(s)} \sup  \vert \gamma_\beta' \vert ^{\Re(s)} $$
$$\times  \left( \int_{\Omega(h)} \vert F(w,\gamma_\alpha g)\vert^2 dm(w) \right)^{1/2} \left( \int_{\Omega(h)} \vert F(w,\gamma_\beta g)\vert^2 dm(w) \right)^{1/2}.$$
Exchanging summations, we can use Cauchy-Schwarz again (and translation invariance of norms with respect to the $g$ variable)  to get
$$\sum_{g \in \G}  \left( \int \vert F(w,\gamma_\alpha g)\vert^2  \right)^{1/2} \left( \int \vert F(w,\gamma_\beta g)\vert^2 \right)^{1/2}
\leq \Vert F \Vert_{H^2_q(h)}^2.$$
This concludes the proof since by Lemma \ref{pressure2}, we now have
$$\Vert \lt_s^N(F)\Vert^2_{H^2_q(h)}\leq C h^{-2\delta}e^{C \vert \Im(s) \vert h} e^{2NP(\Re(s))}.\  \square$$


\begin{thebibliography}{}

\end{thebibliography}


\begin{thebibliography}{10}

\bibitem{Bandtlow1}
Oscar~F. Bandtlow and Oliver Jenkinson.
\newblock On the {R}uelle eigenvalue sequence.
\newblock {\em Ergodic Theory Dynam. Systems}, 28(6):1701--1711, 2008.

\bibitem{Borthwick}
David Borthwick.
\newblock {\em Spectral theory of infinite-area hyperbolic surfaces}, volume
  256 of {\em Progress in Mathematics}.
\newblock Birkh\"auser Boston Inc., Boston, MA, 2007.

\bibitem{Borthwick3}
David Borthwick.
\newblock Sharp geometric upper bounds on resonances for surfaces with
  hyperbolic ends.
\newblock {\em Anal. PDE}, 5(3):513--552, 2012.

\bibitem{BG1}
Jean Bourgain and Alex Gamburd.
\newblock Uniform expansion bounds for {C}ayley graphs of {${\rm SL}_2(\Bbb
  F_p)$}.
\newblock {\em Ann. of Math. (2)}, 167(2):625--642, 2008.

\bibitem{BGS}
Jean Bourgain, Alex Gamburd, and Peter Sarnak.
\newblock Generalization of {S}elberg's {$3/16$} theorem and affine sieve.
\newblock {\em Arxiv preprint}, 2009.

\bibitem{BK1}
Jean Bourgain and Alex Kontorovich.
\newblock On representations of integers in thin subgroups of {${\rm SL}\sb
  2(\bold Z)$}.
\newblock {\em Geom. Funct. Anal.}, 20(5):1144--1174, 2010.

\bibitem{BK2}
Jean Bourgain and Alex Kontorovich.
\newblock On {Z}aremba's conjecture.
\newblock {\em Arxiv preprint}, 2011.

\bibitem{Bowen1}
Rufus Bowen.
\newblock Hausdorff dimension of quasicircles.
\newblock {\em Inst. Hautes \'Etudes Sci. Publ. Math.}, ({\bf 50}):11--25,
  1979.

\bibitem{button}
Jack Button.
\newblock All {F}uchsian {S}chottky groups are classical {S}chottky groups.
\newblock In {\em The Epstein birthday schrift}, volume~{\bf 1} of {\em Geom.
  Topol. Monogr.}, pages 117--125 (electronic). Geom. Topol. Publ., Coventry,
  1998.

\bibitem{Gamburd}
Alex Gamburd.
\newblock On the spectral gap for infinite index ``congruence'' subgroups of
  {${\rm SL}\sb 2(\bold Z)$}.
\newblock {\em Israel J. Math.}, 127:157--200, 2002.

\bibitem{Gohberg}
Israel Gohberg, Seymour Goldberg, and Nahum Krupnik.
\newblock {\em Traces and determinants of linear operators}, volume 116 of {\em
  Operator Theory: Advances and Applications}.
\newblock Birkh\"auser Verlag, Basel, 2000.

\bibitem{GuiZwor2}
L.~Guillop{\'e} and M.~Zworski.
\newblock The wave trace for {R}iemann surfaces.
\newblock {\em Geom. Funct. Anal.}, {\bf 9}(6):1156--1168, 1999.

\bibitem{GLZ}
Laurent Guillop\'e, Kevin~K. Lin, and Maciej Zworski.
\newblock The {S}elberg zeta function for convex co-compact {S}chottky groups.
\newblock {\em Comm. Math. Phys.}, 245(1):149--176, 2004.

\bibitem{GuiZwor1}
Laurent Guillop{\'e} and Maciej Zworski.
\newblock Upper bounds on the number of resonances for non-compact {R}iemann
  surfaces.
\newblock {\em J. Funct. Anal.}, {\bf 129}(2):364--389, 1995.

\bibitem{GuiZwor}
Laurent Guillop{\'e} and Maciej Zworski.
\newblock Scattering asymptotics for {R}iemann surfaces.
\newblock {\em Ann. of Math. (2)}, {\bf 145}(3):597--660, 1997.

\bibitem{JakobsonNaud2}
Dmitry Jakobson and Fr{\'e}d{\'e}ric Naud.
\newblock On the critical line of convex co-compact hyperbolic surfaces.
\newblock {\em Geom. Funct. Anal.}, 22(2):352--368, 2012.

\bibitem{LP2}
Peter~D. Lax and Ralph~S. Phillips.
\newblock Translation representation for automorphic solutions of the
  non-{E}uclidean wave equation {I}, {II}, {III}.
\newblock {\em Comm. Pure. Appl. Math.}, {\bf 37,38}:303--328, 779--813,
  179--208, 1984, 1985.

\bibitem{MazzMel}
Rafe~R. Mazzeo and Richard~B. Melrose.
\newblock Meromorphic extension of the resolvent on complete spaces with
  asymptotically constant negative curvature.
\newblock {\em J. Funct. Anal.}, {\bf 75}(2):260--310, 1987.

\bibitem{Naud2}
Fr{\'e}d{\'e}ric Naud.
\newblock Expanding maps on {C}antor sets and analytic continuation of zeta
  functions.
\newblock {\em Ann. Sci. \'Ecole Norm. Sup. (4)}, 38(1):116--153, 2005.

\bibitem{NaudInvent}
Fr{\'e}d{\'e}ric Naud.
\newblock Density and location of resonances for convex co-compact hyperbolic
  surfaces.
\newblock {\em Invent. Math.}, 195(3):723--750, 2014.

\bibitem{PP}
William Parry and Mark Pollicott.
\newblock Zeta functions and the periodic orbit structure of hyperbolic
  dynamics.
\newblock {\em Ast\'erisque}, (187-188):268, 1990.

\bibitem{Patterson1}
S.~J. Patterson.
\newblock The limit set of a {F}uchsian group.
\newblock {\em Acta Math.}, {\bf 136}(3-4):241--273, 1976.

\bibitem{PatPer}
S.~J. Patterson and Peter~A. Perry.
\newblock The divisor of {S}elberg's zeta function for {K}leinian groups.
\newblock {\em Duke Math. J.}, {\bf 106}(2):321--390, 2001.
\newblock Appendix A by Charles Epstein.

\bibitem{Ruelle0}
David Ruelle.
\newblock Zeta-functions for expanding maps and {A}nosov flows.
\newblock {\em Invent. Math.}, 34(3):231--242, 1976.

\bibitem{Simon}
Barry Simon.
\newblock {\em Trace ideals and their applications}, volume~35 of {\em London
  Mathematical Society Lecture Note Series}.
\newblock Cambridge University Press, Cambridge, 1979.

\bibitem{Tit}
E.~C. Titchmarsh.
\newblock {\em The theory of functions}.
\newblock Oxford University Press, second edition, 1932.

\end{thebibliography}
\end{document}